\documentclass[12pt]{article}

\usepackage{amsmath}
\usepackage{amssymb}
\usepackage{amsthm}
\usepackage{a4wide}
\usepackage[english]{babel}

\numberwithin{equation}{section}

\newtheorem{Th}{Theorem}
\newtheorem{Lem}{Lemma}
\newtheorem{Prop}{Proposition}
\newtheorem{Cor}{Corollary}

\newtheorem{Rem}{Remark}

\newcommand{\Dem}{\noindent{\bf Proof }}

\newcommand{\N}{\mathbb{N}}
\newcommand{\Z}{\mathbb{Z}}
\newcommand{\Q}{\mathbb{Q}}
\newcommand{\Qbar}{\overline{\mathbb{Q}}}
\newcommand{\R}{\mathbb{R}}
\newcommand{\C}{\mathbb{C}}

\newcommand{\etoile}{^\ast}

\newcommand{\Span}{{\rm Span}}

\newcommand{\eps}{\varepsilon}

\newcommand{\om}{\omega}
\newcommand{\eneq}{\end{equation}}
\newcommand{\lam}{\lambda}

\newcommand{\una}{\{1,\ldots,a\}}
\newcommand{\deuxa}{\{2,\ldots,a\}}
\newcommand{\unN}{\{1,\ldots,N\}}
\newcommand{\zeronsN}{\{0,\ldots,n/N\}}

\newcommand{\calN}{\mathcal{N}}

\newcommand{\ord}{{\rm ord}}
\newcommand{\dd}{{\rm d}}
\newcommand{\rk}{{\rm rk}}

\newcommand{\Card}{\#}
\newcommand\tra{ \ ^t  }

\newcommand{\cstnum}{c_1}
\newcommand{\cdeux}{c_2}

\newcommand{\dz}{d_0}

\newcommand{\combitiny}[2]{{\tiny \left( \! \!  \begin{array}{c} #1 \\ #2 \end{array} \! \! \right)}}

\newcommand{\calD}{{\mathcal{D}}}
\newcommand{\calF}{{\mathcal{F}}}

\newcommand{\xii}{L(f,i)}
\newcommand{\xij}{L(f,j)}

\newcommand{\Pbar}{\overline P}

\title{Irrationality of values of $L$-functions of Dirichlet characters}
\author{St\'ephane Fischler}
\date{\today}

\begin{document}

\newcommand{\Li}{{\rm Li}}
\newcommand{\norm}[1]{\lVert #1\rVert}

\newcommand{\ttt}{i_0}
\maketitle

\begin{abstract}
In a recent paper with Sprang and Zudilin, the following result was proved: if $a$ is large enough in terms of $\varepsilon>0$, then at least 
$2^{(1-\varepsilon)\frac{\log a}{\log \log a}}$ values of the Riemann zeta function at odd integers between $3$ and $a$  are irrational. This improves on the Ball-Rivoal theorem, that provides only $\frac{1-\varepsilon}{1+\log 2} \log a$ such irrational values -- but with a stronger property: they are linearly independent over the rationals.

In the present paper we generalize this recent result to both $L$-functions of Dirichlet characters and Hurwitz zeta function. The strategy is different and less elementary: the construction is related to a Pad\'e approximation problem, and a generalization of Shidlovsky's lemma is used to apply Siegel's linear independence criterion. 

We also improve the analogue of the  Ball-Rivoal theorem in this setting: we obtain  $\frac{1-\varepsilon}{1+\log 2} \log a$  linearly independent values $L(f,s)$ with $s\leq a$ of a fixed parity, when $f$ is a Dirichlet character. The new point here is that the constant $1+\log 2$ does not depend on $f$.
\end{abstract}

\noindent{\bf MSC 2010:} 11J72 (Primary); 11M06, 11M35,  33C20 (Secondary).

\bigskip

The purpose of this paper is to prove results of irrationality, or linear independence, of values of the Hurwitz $\zeta$ function or $L$-functions of Dirichlet characters. Both are generalizations of the Riemann $\zeta$ function, so we begin with a quick survey of the main results in this  setting.

When $s\geq 2$ is even, $\zeta(s)\pi^{-s}$ is a non-zero rational number so that $\zeta(s)$ is transcendental. Ap\'ery has proved \cite{Apery} that $\zeta(3)$ is irrational, but there is no odd $s\geq 5$ for which $\zeta(s)$ is known to be irrational. The next breakthrough is due to Ball-Rivoal \cite{BR, RivoalCRAS}:
$$\dim_\Q \Span_\Q (1, \zeta(3), \zeta(5), \ldots, \zeta(a)) \geq \frac{1+o(1)}{1+\log 2}\log a \mbox{ as $a\to\infty$, $a$ odd.}$$
Here and throughout this introduction, $o(1)$ denotes any sequence that tends to 0 as $a\to\infty$. In this paper we mention only asymptotic results (namely, as $a\to\infty$) eventhough most results can be made explicit, and often refined, for small values of $a$. At last we mention the following recent result \cite{FSZ}:
\begin{equation} \label{eqFSZ}
\mbox{at least $2^{(1-o(1))\frac{\log a}{\log \log a}}$ numbers among $ \zeta(3)$,  $ \zeta(5)$,  \ldots, $\zeta(a)$ are irrational,}
 \end{equation}
for $a$ odd, $a\to\infty$.

\bigskip

The natural setting to generalize these results to values of the Hurwitz $\zeta$ function or  $L$-functions of Dirichlet characters is the following. Let $T \geq 1$, and $f: \Z \to \C$ be such that $f(n+T)=f(n)$ for any $n$. We assume that $f$ is not identically zero.  Let  $\eps>0$, and $a $ be sufficiently large (in terms of  $T$ and $\eps$).  For   $p\in\{0,1\}$ consider the complex numbers 
\begin{equation} \label{eq01}
L(f,s) =  \sum_{n=1}^\infty \frac{f(n) }{n^s} \mbox{ with $2\leq s \leq a$ and $s\equiv p \bmod 2$.}
 \end{equation}
If $f$ is a Dirichlet character  mod $T$  then
  these are exactly the values  of the associated $L$-function.

\bigskip

 The restriction on the parity of $s$ in \eqref{eq01} is needed in some cases to get rid of powers of $\pi$. Indeed, if 
 $f$ is a Dirichlet character     then $f$ is either even (i.e., $f(-n)=f(n)$) or odd  (i.e., $f(-n)=-f(n)$), according to whether $f(-1)$ is equal to 1 or $-1$. If $s\geq 2$ has the same parity as $f$ then $L(f,s)\pi^{-s}$ is a non-zero algebraic number (see for instance   \cite[Chapter VII, \S 2]{Neukirch}) so that the numbers  $L(f,s)$ for $s$ with this parity are linearly independent over $\Qbar$. Moreover, for any periodic map $f : \Z \to\Q$ which is either even 
or odd (and not identically zero), we also have  $L(f,s)\pi^{-s}\in\Qbar\etoile$ when $s$ and $f$ have the same parity (see \cite{ChowlaMilnor}). In these situations, we prove new results on the numbers  \eqref{eq01} only when $p$ and $f$ have opposite parities.

\bigskip

An interesting case where (in general) $f$ is neither odd nor even is the following. Given $u\in\{1,\ldots,T-1\}$ we  define $f$ by $f(n)  = 1$ if $n\equiv u \bmod T$, and $f(n)=0$ otherwise.  Then 
 $$L(f,s)  =  \sum_{k=0}^\infty \frac{1 }{(kT+u)^s} =\frac1{T^s}   \sum_{k=0}^\infty \frac{1 }{(k +u/T)^s} 
 = \frac1{T^s}  \zeta(s,\frac{u}{T})$$
 where $\zeta(s,\alpha)$ is  the Hurwitz $\zeta$ function. Therefore the general setting \eqref{eq01} encompasses both values of the Hurwitz $\zeta$ function and values  of $L$-functions of Dirichlet characters.

\bigskip

As far as we know, Ap\'ery's theorem has never been generalized in this direction; the first natural conjecture in this respect is probably that Catalan's constant $L(\chi,2)$ is irrational, where $\chi$ is the non-principal character mod 4. The Ball-Rivoal theorem has been generalized to the  $L$-function of this character  by Rivoal and Zudilin \cite{Catalan}: they have proved \eqref{eqNishi} below with $2+\log 2$ instead of $T+\log 2$, eventhough $T=4$. In the general setting of \eqref{eq01},  Nishimoto has generalized the Ball-Rivoal theorem as follows \cite{Nishimoto}:
\begin{equation} \label{eqNishi}
\dim_\Q \Span_\Q \Big\{ L(f,s), \, \,    2\leq s \leq a, \, \, s\equiv p \bmod 2 \Big\} 
  \geq \frac{1+o(1)}{T+\log 2}\log a \mbox{ as $a\to\infty$.}
 \end{equation}
In the special case where $\sum_{n=1}^T f(n) \neq 0$ (which includes the Hurwitz $\zeta$ function but not $L$-functions of non-principal Dirichlet characters), this lower bound appears already in Nash' thesis \cite{Nash}. The constant $ T+\log 2 $ in Eq. \eqref{eqNishi} has been refined to  $ T/2 +\log 2 $ in \cite{SFcaract}, provided $f$ is a Dirichlet character and $T$ is a multiple of 4. When $f$ is the non-principal character mod 4, this gives as a special case the lower bound of Rivoal and Zudilin \cite{Catalan}.

\bigskip

Our first result is that one may replace the constant  $ T+\log 2 $ in Eq. \eqref{eqNishi}  with   $ 1+\log 2 $, so that the lower bound is uniform in $T$  and  is the same as for the Riemann $\zeta$ function.

\begin{Th} \label{th1}
Let $T \geq 1$, and $f: \Z \to \C$ be such that $f(n+T)=f(n)$ for any $n$. Assume that $f$ is not identically zero.  Let $p\in\{0,1\}$, $\eps>0$, and $a $ be sufficiently large (in terms of    $T$  and $\eps$).  Then 
$$\dim_\Q \Span_\Q \Big\{ L(f,s), \, \,    2\leq s \leq a, \, \, s\equiv p \bmod 2 \Big\} 
  \geq \frac{1-\eps}{1+\log 2}\log a.
  $$
\end{Th}

Of course the same result holds  without the  restriction $s\equiv p \bmod 2$, but it is weaker and even trivial in some cases where $f$ is even or odd (as noticed above). 

\bigskip

In another direction,  we generalize the recent result \eqref{eqFSZ} to this setting.

\begin{Th} \label{th2}
Let $T \geq 1$, and $f: \Z \to \C$ be such that $f(n+T)=f(n)$ for any $n$. Assume that $f$ is not identically zero.  Let $E$ be a finite-dimensional $\Q$-vector space contained in $\C$, $p\in\{0,1\}$, $\eps>0$, and $a $ be sufficiently large (in terms of $\dim E$, $T$, and $\eps$).  Then among the numbers
$  L(f,s) $  with $2\leq s \leq a$ and $s\equiv p \bmod 2$,
at least
$$2^{(1-\eps)\frac{\log a}{\log \log a}}$$
do not belong to $E$.
\end{Th}

Taking $E=\Q$ we obtain at least $2^{(1-\eps)\frac{\log a}{\log \log a}}$ irrational values among the numbers $L(f,s)$. The dependence in $a$ is much better than in the lower bound of Theorem \ref{th1}; however we obtain only numbers outside $E$, and not $\Q$-linearly independent numbers.

\bigskip

Before explaining the strategy used in the proofs of Theorems \ref{th1} and \ref{th2}, we would like to state the two main special cases of Theorem  \ref{th2} explicitly. 

\begin{Cor} \label{cor1}
Let $\chi$ be a Dirichlet character; put $p=0$ is $\chi $ is odd, and $p=1$ if $\chi$ is even.  Let $E$ be a finite-dimensional $\Q$-vector space contained in $\C$. Let   $\eps>0$, and $a $ be sufficiently large (in terms of $\chi$, $\dim E$,  and $\eps$).  Then among the numbers $L(\chi,s)$ with   $2\leq s \leq a$ and $s\equiv p \bmod 2$,
at least
 $2^{(1-\eps)\frac{\log a}{\log \log a}}$ 
do not belong to $E$.
\end{Cor}

\begin{Cor} \label{cor2} Let $r$ be a positive rational number, and $p\in\{0,1\}$. 
 Let $E$ be a finite-dimensional $\Q$-vector space contained in $\C$. Let   $\eps>0$, and $a $ be sufficiently large (in terms of $r$, $\dim E$,  and $\eps$).  Then among the numbers 
 $$\zeta(s,r) = \sum_{n=0}^\infty \frac1{(n+r)^s}$$
  with   $2\leq s \leq a$ and $s\equiv p \bmod 2$,  at least
 $2^{(1-\eps)\frac{\log a}{\log \log a}}$ 
do not belong to $E$.
\end{Cor}

Corollary \ref{cor2} is new even for $r=1$, i.e. for the Riemann $\zeta$ function: it is a refinement of \eqref{eqFSZ}. We would like to emphasize the fact that the proof of \cite{FSZ} {\em does not} give this result for $E \neq \Q$: a different approach is used here, proving linear independence and not only irrationality.

\bigskip

The proof  of Theorems \ref{th1} and \ref{th2} is based on the strategy of \cite{SFcaract}: we apply Siegel's linear independence criterion using a general version of Shidlovsky's lemma (namely   Theorem~\ref{thshid}, stated in \S \ref{subsecshid} and proved in \cite{SFcaract} following the approach  of Bertrand-Beukers \cite{BB} and Bertrand \cite{DBShid}). This makes it necessary to relate the construction to a Pad\'e approximation problem with essentially as many equations as the number of unknowns. In the present paper we adapt this strategy so as to include Sprang's arithmetic lemma \cite[Lemma 1.4]{Sprang} and the elimination trick of \cite{Zudilintrick, Sprang, FSZ}. The proofs  of Theorems \ref{th1} and \ref{th2} are essentially the same, except for the choice of parameters. It is also  possible to prove other results of the same flavour  (see Theorem \ref{th3} at the end of \S \ref{subsecdemth2}, which implies both Theorem \ref{th2}  and -- up to a multiplicative constant -- Theorem \ref{th1}). 

Our construction contains as a special case the one used in \cite{FSZ} to prove \eqref{eqFSZ}. We prove this in \S \ref{subsecFSZ}; as a byproduct, we relate the construction of \cite{FSZ} to a Pad\'e approximation problem with essentially as many equations as the number of unknowns. 

\bigskip

The structure of this paper is as follows. We gather in Proposition \ref{propdio} the output of the Diophantine construction (see \S \ref{subsec21}), and prove it in \S \ref{secdio}. Then we deduce Theorems \ref{th1} and~\ref{th2} from Proposition~\ref{propdio} in  \S \ref{sec1}  using Siegel's linear independence criterion (stated in  \S \ref{sec3}).

\section{Diophantine construction} \label{secdio}

In this section we gather the  Diophantine part of the proof, namely the  construction of linearly independent linear forms. We prove
   Proposition~\ref{propdio} stated in \S \ref{subsec21}, from which we  shall  deduce  in  \S \ref{sec1} the results stated in the introduction. The linear forms are constructed in \S \ref{subsec33} using series of hypergeometric type. We relate them in \S \ref{subsec34} to a Pad\'e approximation problem, and then apply a general version of Shidlovsky's lemma (stated in \S \ref{subsecshid}). At last, arithmetic and asymptotic properties are dealt with in \S \ref{subsec35}.

\subsection{Statement of the result}  \label{subsec21}

Let $a$, $r$, $N$ be positive integers such that $1 \leq r < \frac{a}{3N}$. 
Let $N \geq 1$, and $f: \Z \to \C$ be such that $f(m+N)=f(m)$ for any $m$. Assume that $f$ is not identically zero. Let $p\in\{0,1\}$; put
$$\xij = \sum_{m=1}^\infty \frac{f(m) }{m^j} \mbox{ for any } j \in\deuxa.$$
Let also
\begin{equation} \label{eqdefalphabeta}
\alpha = (4e)^{(a+1)/N} (2N)^{2r+2} r^{-(a+1)/N+4(r+1)} \mbox{ and }
\beta = (2e )^{(a+1)/N} (r +1)^{2r+2} N^{2r+2}.
\eneq

\begin{Prop} \label{propdio}
There exists a constant $ \cstnum$, which depends only on $a$ and $N$, with the following property. For any integer multiple $n$ of $N$ there exist integers $s_{k,i}$, with $1\leq k \leq  \cstnum$ and $2\leq i \leq a+N$, such that:
\begin{enumerate}
\item[$(i)$] For any $n$ sufficiently large, the subspace $\calF$ of $\R^{a+N-1}$ spanned by the vectors \linebreak $\tra (s_{k,2},\ldots,s_{k,a+N})$, $1\leq k \leq \cstnum$, is non-zero and does not depend on $n$.
\item[$(ii)$] For any $k$ and any $i$ we have $|s_{k,i}| \leq \beta^{n+o(n)}$ as $n\to\infty$.
\item[$(iii)$] For any $k$  we have, as $n\to\infty$:
\begin{equation} \label{eqnvprop}
\Big| 2(-1)^p \sum_{\substack{i=2 \\ i\equiv p \bmod 2}}^a s_{k,i} \xii + \sum_{i=0}^{N-1} s_{k,a+1+i}f(i)\Big| \leq \alpha ^{n+o(n)}.
\eneq
\end{enumerate}
\end{Prop}

From now on, the symbols $o(\cdot)$ will be intended as $n\to\infty$.  
 Since $k \leq \cstnum$, these symbols   can be made uniform with respect to $k$.

  \bigskip
  
  The integers $s_{k,i}$ depend also implicitly on $n$, $a$, $r$, $N$, $f$ and $p$. Their values for $i \not\equiv p \bmod 2$ do not appear in the linear combinations of part $(iii)$, but they could be of interest in other settings. Another feature of this construction is that for $i\leq a$,  the integers $s_{k,i}$  depend  only on $n$, $a$, $r$, $N$ but not on $f$ or $p$. Probably this could lead to variants of our results in the style of \cite{Pilehroodssums} or \cite{SFdistrib}.

  \bigskip

\begin{Rem} \label{remenonce} In \cite{SFcaract} a similar construction is made, where the  matrix $[s_{k,i}]_{i,k}$ has rank $a+N-1$ for $n$ sufficiently large so that the subspace $\calF$ of part $(i)$ is equal to $\R^{a+N-1}$. In the present setting we make a different construction to incorporate Sprang's arithmetic lemma (see \S \ref{subsec33} below), and the matrix  $[s_{k,i}] $ we obtain has rank less than  $a+N-1$  for some values of the parameters (see Remark \ref{remf} in \S \ref{subsec34}): the  subspace $\calF$  in Proposition \ref{propdio} is not always equal to $\R^{a+N-1}$. 
\end{Rem}

\subsection{Construction of the linear forms}\label{subsec33}

In this section we define the numbers $s_{k,i}$ of Proposition \ref{propdio} (see Eq. \eqref{eqdefski}) and express the linear form of Eq. \eqref{eqnvprop} in a more convenient way (see Lemma  \ref{lemdcp}).  We postpone until \S \ref{subsec35} the proof that $s_{k,i}\in\Z$. 

\bigskip 

As in \S \ref{subsec21} we let $a$, $r$, $N$ be positive integers such that $1 \leq r < \frac{a}{3N}$. For any integer multiple $n$ of $N$ we let
$$F(t) = (n/N)!^{(a+1)-(2r+1)N} \frac{(t-rn)_{(2r+1)n+1}  }{\prod_{h=0}^{n/N}(t+Nh)^{a+1}}$$
where $(\alpha)_p = \alpha (\alpha+1)\ldots (\alpha+p-1)$ is the Pochhammer symbol. 
Note that each factor $t+Nh$ of the denominator appears also in the numerator, so that the poles $t = -Nh$ of this rational function only have  order $a$.  This rational function $F(t)$ is  similar to that of \cite{SFcaract}, but central factors have been inserted in the numerator to apply Sprang's arithmetic lemma (see Remark \ref{remnv} below). 

\bigskip

In this section we follow the proof of \cite{SFcaract}, except for Eq. \eqref{equplusv} which is specific to the function $F$ we consider here. We let 
\begin{equation} \label{eqdefszinf}
S_0(z) = \sum_{t=n+1}^\infty F(-t) z^t \hspace{1cm}  \mbox{ and }  \hspace{1cm}  
S_\infty(z) = \sum_{t=1}^\infty F( t) z^{-t} 
\eneq
for $z\in\C$ with $|z|=1$; then both series are convergent since  the degree $-\dz$ of   $F$   satisfies 
\begin{equation} \label{eqdefdz}
\dz :=  -\deg F = (a+1)  (\frac{n}{N} +1) - (2r+1)n -1   \geq 2.
\eneq
We let $\om = e^{2i\pi/N}$ and for any $\ell\in\unN$ we consider the (inverse) discrete Fourier transform of $f$, defined by 
\begin{equation} \label{eqdefmunu}
\widehat{f}(\ell) = \frac1{N} \sum_{\lam = 1}^N f(\lam) \om^{-\ell \lam}.
\eneq
We also let 
$$\delta_n = (N d_{n/N})^{a+1} N^{(a+1)n/N}, \mbox{ where } d_j = {\rm lcm}(1,2,\ldots,j).$$
The linear forms of Proposition \ref{propdio} are given by the following lemma. The rational numbers $s_{k,i}$ will be constructed explicitly in the proof (see Eq. \eqref{eqdefski}), and we shall prove in \S \ref{subsec35}   that  they are integers.

\begin{Lem}\label{lemdcp} For any $1\leq k \leq \dz-1$ there exist rational numbers $s_{k,2}$, \ldots, $s_{k, a+N}$ such that 
\begin{eqnarray}
&&
 \delta_n  \sum_{\ell=1}^N \widehat{f}(\ell)  \Big[  \om^{\ell (k-1)} S_0^{(k-1)}(\om^\ell  ) + (-1)^p   \om^{\ell (1-k)}S_\infty^{(k-1)}( \om^{-\ell
} )\Big] \nonumber 
\\ 
&&= 
2(-1)^p \sum_{\substack{i=2 \\ i\equiv p \bmod 2}}^a s_{k,i} \xii + \sum_{i=0}^{N-1} s_{k,a+1+i}f(i) 
 \label{eqlemdcp}
\end{eqnarray}
where   $S^{(k-1)}$ is the $(k-1)$-th derivative of a function $S$.
 \end{Lem}

\bigskip

Let us prove Lemma \ref{lemdcp}. The partial fraction expansion of $F$ reads
$$F(t) = \sum_{h=0}^{n/N}\sum_{j=1}^a \frac{p_{j,h}}{(t+Nh)^j}$$
with rational coefficients $p_{j,h}$; we consider 
$$P_j(z) = \sum_{h=0}^{n/N} p_{j,h} z^{Nh} \in\Q[z]_{\leq n} \mbox{ for any } j\in\{1,\ldots,a\}.$$
Let   $P_{1,j} = P_j$ for any $j\in\una$, and define inductively $P_{k,j}\in\Q(z)$ by
\begin{equation} \label{eqdefpkj}
P_{k,j}(z) = P'_{k-1,j}(z) - \frac1{z} P_{k-1,j+1}(z) \mbox{ for any } k\geq 2 \mbox{ and any } j\in\una,
\eneq
where   $P_{k-1,a+1} = 0$ for any $k$. We let also\footnote{There is a misprint in the formula  that gives $U(z)$ in \cite{SFcaract}.}
\begin{equation} \label{eqdefu1}
U_1(z) =   -\sum_{t=1}^n z^t \sum_{j=1}^a \sum_{h=0}^{\lfloor (t-1)/N\rfloor} \frac{p_{j,h}}{(Nh-t)^j} \in\Q[z]_{\leq n}
\eneq
\begin{equation} \label{eqdefv1}
\mbox{ and } V_1(z) =    -\sum_{t=0}^{n-1} z^t \sum_{j=1}^a \sum_{h= \lceil (t+1)/N\rceil} ^{n/N}  \frac{p_{j,h}}{(Nh-t)^j} \in\Q[z]_{\leq n},
\eneq
and define $U_k$, $V_k$ for any $k\geq 2$ by the recurrence relations
\begin{equation} \label{eqdefpz}
U_k(z) = U'_{k-1}(z) - \frac1{1-z} P_{k-1,1}(z),
\eneq
\begin{equation} \label{eqdefpinf}
V_k(z) = V'_{k-1}(z) + \frac1{z(1-z)} P_{k-1,1}(z).
\eneq
Then for any $k\geq 1$ we have (as in \cite{BR, FR})
\begin{equation} \label{eqszeroder}
S_0^{(k-1)}(z) = U_k(z) + \sum_{j=1}^a P_{k,j}(z) (-1)^j \Li_j( z)
\eneq
\begin{equation} \label{eqsinfder}
\mbox{ and } S_\infty^{(k-1)}(z) = V_k(z) + \sum_{j=1}^a P_{k,j}(z)  \Li_j(1/z).
\eneq

Since $P_j(z)\in\Q[z^N]$ for any $j\in\una$, Eq. \eqref{eqdefpkj} yields $P_{k,j} \in z^{1-k} \Q[z^N]$. This property is  very important to us  since we shall evaluate $P_{k,j}$ at $N$-th roots of unity. To evaluate in the same way the rational functions  
 $U_k,V_k \in \Q[z,z^{-1}]$ for   $k\leq \dz-1$, we  write
\begin{equation} \label{eqpkzinfom}
z^{ k-1 } U_k(z) = \sum_{\lam = 0}^{N-1} z^\lam U_{k,\lam}(z) 
\mbox{ and }
z^{ k-1 }V_k(z) = \sum_{\lam = 0}^{N-1} z^\lam V_{k,\lam}(z) 
\eneq
with $U_{k,\lam},  V_{k,\lam} \in  \Q[z^N,z^{-N}]$.  Then  Eqns. \eqref{eqszeroder} and  \eqref{eqsinfder} yield
\begin{equation} \label{eqszeroderbis}
z^{k-1} S_0^{(k-1)}(z) =\sum_{\lam = 0}^{N-1} z^\lam U_{k,\lam}(z)  + \sum_{j=1}^a z^{k-1}  P_{k,j}(z) (-1)^j \Li_j( z)
\eneq
\begin{equation} \label{eqsinfderbis}
\mbox{ and } z^{k-1}  S_\infty^{(k-1)}(z) = \sum_{\lam = 0}^{N-1} z^\lam V_{k,\lam}(z)  + \sum_{j=1}^a z^{k-1}  P_{k,j}(z)  \Li_j(1/z).
\eneq

We may  now define the coefficients $s_{k,i}$ for any $k\geq 1$ by:
\begin{equation} \label{eqdefski}
\left\{
\begin{array}{l}
s_{k,i} =  \delta_n  P_{k,i}(1) \mbox{ for } 2\leq i \leq a, \\
s_{k,a+1+\lam} = \delta_n  (U_{k,\lam}(1)+(-1)^p V_{k,N-\lam}(1))  \mbox{ for } 0\leq \lam \leq N-1,
\end{array}
\right.
\eneq
where $V_{k,N} = V_{k,0}$; recall that 
$\delta_n = (N d_{n/N})^{a+1} N^{(a+1)n/N} $ with $d_j = {\rm lcm}(1,2,\ldots,j)$. Since  $ P_{k,j}(z)$, $U_{k,\lam}(z)$ and $V_{k,N-\lam}(z)$ are polynomials  with rational coefficients, the  numbers $s_{k,2}$, \ldots, $s_{k, a+N}$ are rational. We shall prove in Lemma \ref{lemtech}   that they are integers,  thanks to the factor $\delta_n$. We also point out that $s_{k,i}$ is not defined for $i=1$; actually $P_{k,1}(1)=0$ for the values of $k$ we are interested in (see \eqref{eqokenun} below). 

\bigskip

To conclude the proof of Lemma \ref{lemdcp}, we shall evaluate  Eqns. \eqref{eqszeroderbis} and \eqref{eqsinfderbis} at roots of unity. At the point 1 this is possible since, as  in \cite[\S 4.3]{SFcaract},
\begin{equation} \label{eqokenun}
\mbox{for any }k\leq \dz-1,  \hspace{.7cm} P_{k,1}(1) = 0  \hspace{.2cm} \mbox{ and }  \hspace{.2cm} U_k, V_k \mbox{ do not have a pole at } z=1.
\eneq
Now let $k\leq \dz-1$, and $z\in\C$ be such that $|z|=1$. Then Eqns. \eqref{eqszeroder} to  \eqref{eqsinfderbis} hold, upon agreeing that the sums start at $j=2$ if $z=1$; this  remark will be used below when $z$ is a $N$-th root of unity.

Let $\Lambda_k$ be the right hand side of Eq. \eqref{eqlemdcp}. Using  \eqref{eqokenun} the definition \eqref{eqdefski} of $s_{k,i}$ yields 
$$\Lambda_k = 2 \delta_n(-1)^p \sum_{1\leq j \leq a \atop j\equiv p \bmod 2} P_{k,j}(1) \xij +  \delta_n\sum_{\lam=0}^{N-1} (U_{k,\lam}(1)+(-1)^p V_{k,N-\lam}(1))  f(\lam).$$
Now Eq. \eqref{eqdefmunu} yields
$$  \sum_{\ell=1}^N \widehat{f}(\ell) \om^{m \ell} = f(m) \mbox{ for any $m\in\Z$, so that }
  \sum_{\ell=1}^N \widehat{f}(\ell) \Li_j(\om^\ell ) = \sum_{m=1}^\infty \frac{f(m) }{m^j} =  \xij \mbox{ for any } j \leq a.$$ 
  Therefore we have, since $V_{k,N} = V_{k,0}$:
\begin{eqnarray*}
\Lambda_k 
=  
 &&  
  \delta_n \sum_{j=1}^a    P_{k,j}(1)   ( (-1)^j + (-1)^p)  \sum_{\ell=1}^N  \widehat{f}(\ell) 
 \Li_j( \om^\ell )\\
 &&
+
 \delta_n\sum_{\lam = 0}^{N-1}  \Big[ \Big(  \sum_{\ell=1}^N  \widehat{f}(\ell) \om^{\ell   \lambda}\Big) 
U_{k,\lam}(1) + (-1)^p    \Big(  \sum_{\ell=1}^N \widehat{f}(\ell)  \om^{- \ell   \lambda}\Big) V_{k,\lam}(1) \Big] 
.
\end{eqnarray*}
Then Eqns.  \eqref{eqszeroderbis} and  \eqref{eqsinfderbis} yield,  since $U_{k,\lam}(z)  $, $V_{k,\lam}(z)  $, and $z^{k-1}  P_{k,j}(z)$ depend only on $z^N$ and $\omega$ is a $N$-th root of unity:
$$\Lambda_k = \delta_n \sum_{\ell=1}^N\widehat{f}(\ell) \Big[  \om^{\ell (k-1)} S_0^{(k-1)}(\om^\ell  ) + (-1)^p   \om^{\ell (1-k)}S_\infty^{(k-1)}( \om^{-\ell
} )\Big].$$
This concludes the proof of Lemma \ref{lemdcp}.

\bigskip

\begin{Rem} \label{remnv} The only difference here with the construction of \cite{SFcaract} is that the rational function $F$ has been modified to incorporate Sprang's arithmetic lemma \cite[Lemma 1.4]{Sprang}. In our setting this choice of $F$ leads to the following additional property, that will be used in \S \ref{subsec34}:
\begin{equation} \label{equplusv}
U_1(z) + V_1(z) \in \Q[z^N].
\eneq
To prove this property  we notice that 
$$U_1(z) + V_1(z) = - \sum_{t=1}^n z^t \sum_{j=1}^a \sum_{h \neq t/N} \frac{p_{j,h}}{(Nh-t)^j} ; $$
for any $t$ which is not a multiple of $N$, the coefficient of $z^t$ is $-F(-t) = 0$. 
\end{Rem}

\subsection{A general version of Shidlosvky's lemma} \label{subsecshid}

Let $q$ be a positive integer, and $A \in M_q(\C(z))$. We fix\footnote{We shall check in \S \ref{subsec34} that the notation introduced in the present section is consistent with the one used earlier in this paper.}  $P_1,\ldots,P_q\in\C[z]$ and $n \in\N=\{0,1,2,\ldots\}$ such that 
  $\deg P_i \leq n$ for any $i$.  Then with any solution $Y = \tra (y_1,\ldots,y_q)$ of the differential system $Y'=AY$   is  associated a remainder $R(Y)$ defined   by 
$$
R(Y)(z)  = \sum_{i=1}^q P_i(z) y_i(z).$$
Let $\Sigma$ be a finite subset of $  \C\cup\{\infty\}$, which may contain singularities of   the differential system $Y'=AY$.
 For each $\sigma\in\Sigma$, let 
$(Y_j)_{j\in J_\sigma}$ be a family of solutions of $Y'=AY$ such that the functions $R(Y_j)$, $j\in J_\sigma$, are $\C$-linearly independent and  
belong to the Nilsson class at $\sigma$ (i.e., have a local expression at $\sigma$ as linear combination of holomorphic functions, with coefficients involving complex powers of $z-\sigma$ and integer powers of $\log(z-\sigma)$). 
We agree that $J_\sigma=\emptyset$ if  $\sigma\not\in    \Sigma$, and define rational functions    $P_{k,i}\in\C(z)$ for $k \geq 1$ and $1 \leq i \leq q$  by
\begin{equation} \label{eqdefpkishid}
\left( \begin{array}{c} P_{k,1} \\ \vdots \\ P_{k,q} \end{array}\right) = \left(\frac{\dd}{\dd z} + \tra A\right)^{k-1} 
 \left( \begin{array}{c} P_{1} \\ \vdots \\ P_{q} \end{array}\right).
 \eneq
  These rational functions $P_{k,i}$ play an important role because they are used to differentiate the remainders  (see   \cite[Chapter 3, \S 4]{Shidlovski}):
\begin{equation} \label{eqderiR}
R(Y)^{(k-1)}(z) =  \sum_{i=1}^q P_{k,i}(z) y_i(z). 
\end{equation}
The following result is proved in \cite[Theorem 1.2]{SFcaract}.

\begin{Th} \label{thshid}
There exists  a  positive constant   $\cdeux$, which depends only on    $A$ and $\Sigma$, with the following property. Assume that, for some $\alpha\in\C$:
  \begin{itemize}
\item[$(i)$] The differential system $Y'=AY$ 
 has a basis of local solutions at $\alpha$ with coordinates in \mbox{$\C[\log(z-\alpha)][[(z-\alpha)^e]]$} for some positive real  number $e$.
 \item[$(ii)$]  We have 
$$
\sum_{\sigma\in\Sigma}\sum_{j\in J_{\sigma}} \ord_\sigma(R(Y_j)) \geq (n+1)  q - n  \Card J_\infty  -\tau 
$$
 for some $\tau $ with $0 \leq \tau \leq n - \cdeux$.
\item[$(iii)$] All rational functions  $P_{k,i}$, with $1\leq i \leq q$ and $1 \leq k < \tau + \cdeux$, are holomorphic at $z=\alpha$.
\end{itemize}
Then the matrix  $[P_{k,i}(\alpha)]_{1\leq i  \leq q, 1\leq k < \tau +  \cdeux} \in M_{q,\tau + \cdeux-1}(\C)$ has rank at least $q-\Card J_\alpha$.
\end{Th}
In the special case where $\Sigma = \{0\}$, $\Card J_0=1$,  $Y_j$ is analytic at $0$, and $\alpha\neq 0$ is not a singularity of the  differential system $Y'=AY$,  this result was used by  Shidlovsky   to prove the Siegel-Shidlovsky theorem on values of $E$-functions (see \cite[Chapter 3, Lemma~10]{Shidlovski}). The functional part of Shidlovsky's lemma has been generalized by  Bertand-Beukers \cite{BB}  to the case where  $\Sigma\subset\C$, $\Card J_\sigma=1$ for any $\sigma$, and all functions $Y_j$ are obtained by analytic continuation from a single one, analytic at all $\sigma\in\Sigma$. Then Bertrand has allowed  \cite[Th\'eor\`eme 2]{DBShid} an arbitrary number of solutions at each $\sigma$,  assuming  that $\infty\not\in\Sigma$ and  the functions $Y_j$, $j \in J_\sigma$, are analytic at  $\sigma$. The proof  \cite{SFcaract} of Theorem \ref{thshid} follows the approach of Bertand-Beukers  and Bertrand, based on differential Galois theory. 

An important feature of Theorem \ref{thshid} is that $\alpha$ may be  a singularity of the differential system $Y'=AY$, and/or an element of   $\Sigma$. Both happen in the present paper, where $\alpha=1$ (see \S \ref{subsec34} where Theorem \ref{thshid} is applied to prove Lemma \ref{lemmatinv}). 
 If $\alpha\not\in\Sigma$ then $J_\alpha=\emptyset$ so that Theorem \ref{thshid} yields  a matrix of maximal rank $q$. On the other hand, if $\alpha\in\Sigma$  then the $\Card J_\alpha$ linearly independent linear combinations of the rows of the matrix $[P_{k,i}(z)]_{i,k}$  corresponding to $R(Y_j)$, $j\in J_\alpha$, may  vanish at $\alpha$: the lower bound $q-\Card J_\alpha$ is best possible. In the setting of  \S \ref{subsec34}  we have $\alpha=1$  and $J_1 = \{1\}$  so that  Theorem \ref{thshid} yields $\rk[P_{k,i}(1)]\geq q-1$. Now \eqref{eqokenun} in the proof of Lemma \ref{lemdcp} shows that $P_{k,1}(1)=0$ for any $k< \tau + \cdeux$ (since $\tau+\cdeux\leq \dz$ because $\tau$ and $\cdeux$ are independent from $n$ whereas $\dz$ tends to $\infty$ with $n$). Therefore the matrix $[P_{k,i}(1)]_{1\leq i \leq q, 1\leq k < \tau + \cdeux}$ has rank equal to $q-1$. Removing the first row, which is identically zero, yields a matrix of rank $q-1$ equal to the number of rows.

\subsection{Pad\'e approximation and application of Shidlovsky's lemma}\label{subsec34}

In this section we prove part $(i)$ of Proposition~\ref{propdio} for the numbers  $s_{k,i}$ constructed in \S \ref{subsec33}. 

\begin{Lem} \label{lemmatinv} Let $s_{k,i}$ be defined by Eq. \eqref{eqdefski}. Then there exists  a positive constant $\cstnum$ (which depends only on $a$ and $N$)  such that 
for any $n$ sufficiently large, the subspace $\calF$ of $\R^{a+N-1}$ spanned by the vectors $\tra (s_{k,2},\ldots,s_{k,a+N})$, $1\leq k \leq \cstnum$, is non-zero and does not depend on $n$.
\end{Lem}

The proof of Lemma \ref{lemmatinv} falls into 3 steps. To begin with,  we construct a Pad\'e approximation problem related to our construction, with essentially as many equations as the number of unknowns; notice that this problem is not the same as in   \cite{SFcaract}, since the function $F$ used in the construction  is   different. Then we  apply a general version of Shidlovsky's lemma, namely Theorem \ref{thshid} stated in \S \ref{subsecshid}. This provides a matrix $P$ with linearly independent rows. At last, we relate the numbers $s_{k,i}$ to $P$ by constructing a matrix $M$ such that $[s_{k,i}]_{i,k} = MP$. The point is that $M$ does not depend on $n$ (whereas $P$ and $[s_{k,i}]$ do). The subspace spanned by the columns of $[s_{k,i}]_{i,k} $ is   the same as the one  spanned by the columns of $M$: it does not depend on $n$.

\bigskip

\noindent {\bf Step 1:} Construction of the Pad\'e approximation problem.

\bigskip

We  recall from \S \ref{subsec33} that 
$$F(t) = (n/N)!^{(a+1)-(2r+1)N} \frac{(t-rn)_{(2r+1)n+1}  }{\prod_{h=0}^{n/N}(t+Nh)^{a+1}}, $$
$$S_0(z) = \sum_{t=n+1}^\infty F(-t) z^t =  U_1(z) + \sum_{j=1}^a P_j(z) (-1)^j \Li_j(z), $$
$$\mbox{ and } S_\infty(z) = \sum_{t=1}^\infty F( t) z^{-t} =  V_1(z) + \sum_{j=1}^a P_j(z) \Li_j(1/z).$$
Since $P_j(z)  \in \C[z^N]$ for any $j\in\una$, we have $P_j(\om^\ell z) = P_j(z)$ for any $\ell\in\Z$. Therefore letting
\begin{equation} \label{eqc111}
R_{0,\ell}(z) = S_0(\om^\ell z),
\hspace{0.6cm} 
R_{\infty,\ell}(z) = S_\infty(\om^\ell z),
\hspace{.6cm} 
\Pbar_{0,\ell}(z) = U_1(\om^\ell z),
\hspace{.6cm} 
\Pbar_{\infty,\ell}(z) = V_1(\om^\ell z) 
\eneq
for any $\ell\in\unN$, we have
\begin{equation} \label{eqP1}
R_{0,\ell}(z) = \Pbar_{0,\ell}(z) +   \sum_{j=1}^a P_j(z) (-1)^j \Li_j(\om^\ell z) = O(z^{(r+1)n+1}), \hspace{0.5cm} z\to 0,
\eneq
\begin{equation} \label{eqP2}
\mbox{ and }
R_{\infty,\ell}(z) = \Pbar_{\infty,\ell}(z) +   \sum_{j=1}^a P_j(z)   \Li_j(\frac1{\om^\ell z}) = O(z^{ -rn-1}), \hspace{0.5cm} z\to \infty.
\eneq
Moreover, recall that $\dz = -\deg F =  (a+1)  (\frac{n}{N} +1) - (2r+1)n -1 $;  Lemma 3 of \cite{FR} shows that
$$
  \sum_{j=1}^a P_j(z)  (-1)^{j-1}\frac{(\log z)^{j-1}}{(j-1)!} = O((z-1)^{\dz-1}),  \hspace{0.5cm} z\to 1.$$
Using again the fact that $  P_j(\om^{-\ell}  z) = P_j(z)$, we obtain for any $\ell \in\unN$:
\begin{equation} \label{eqP3}
R_{\om^\ell  }(z) :=    \sum_{j=1}^a P_j(z)  (-1)^{j-1}\frac{(\log (\om^{-\ell}z))^{j-1}}{(j-1)!} = O((z-\om^\ell)^{\dz-1}),  \hspace{0.5cm} z\to \om^\ell.
\eneq

The new point here, with respect to \cite{SFcaract}, is that Eq. \eqref{equplusv} in Remark \ref{remnv} shows that $\Pbar = \Pbar_{0,\ell} + \Pbar_{\infty,\ell}$ does not depend on $\ell$. Therefore Eq. \eqref{eqP1} can be written as
\begin{equation} \label{eqP1bis}
R_{0,\ell}(z) = \Pbar(z) -   \Pbar_{\infty,\ell}(z) +   \sum_{j=1}^a P_j(z) (-1)^j \Li_j(\om^\ell z) = O(z^{(r+1)n+1}), \hspace{0.5cm} z\to 0.
\eneq

\bigskip

We have obtained a Pad\'e approximation problem with $(a+N+1)(n+1)$ unknowns, namely the coefficients of $\Pbar(z)$, $P_j(z)$ for $1\leq j \leq a$, and $\Pbar_{\infty,\ell}$ for $1\leq \ell \leq N$. Eqns. \eqref{eqP2}, \eqref{eqP3} and \eqref{eqP1bis} amount to 
$$N(( r+1)n+1) + N( \dz-1) + N((r+1)n+1) = (a+N+1)(n+1)-\tau$$
linear equations, where $\tau = a+1-aN$ is the 
 difference between the number of unknowns and the number of equations. If $N=1$ then $\tau=1$: this is exactly the Pad\'e approximation problem of \cite[Theorem 1]{FR}, which has a unique solution up to proportionality. However if $N\geq 2$ then $\tau<0$: we have solved a linear system with (slightly) more equations than the number of unknowns.

\bigskip

\noindent {\bf Step 2:} Application of  Shidlovsky's lemma.

\bigskip

Let us introduce some notation to fit into the context of \S \ref{subsecshid}, and check the assumptions of Theorem \ref{thshid}. 
Let $q = a+N+1$, and $A\in M_q(\C(z))$ be the matrix of which the coefficients $A_{i,j}$ are given by:
$$
\left\{
\begin{array}{l}
A_{i,i-1}(z) = \frac{-1}{z} \mbox{ for any } i \in \{2,\ldots,a\}\\
A_{1,a+ 1}(z) =\frac{1}{z} \\
A_{1,a+1+\ell}(z) = \frac{1}{ z(1-\om^\ell z  )  } \mbox{ for any } \ell\in\unN
\end{array}
\right.
$$
and all other coefficients are zero. We consider the following solutions of the differential system $Y'=AY$, with $1\leq\ell\leq N$:
$$Y_{0,\ell}(z) = \tra \Big( -\Li_1(\om^\ell z), \Li_2(\om^\ell z), \ldots, (-1)^a \Li_a(\om^\ell z), 1, 0, \ldots, 0,-1,0, \ldots,0\Big),$$
$$Y_{\infty,\ell}(z) = \tra \Big(  \Li_1(\frac1{\om^\ell z}), \Li_2(\frac1{\om^\ell z}), \ldots,   \Li_a(\frac1{\om^\ell z}),0, 0,  \ldots, 0,1,0, \ldots,0\Big),$$
$$Y_{\om^\ell}(z) = \tra \Big( 1, -\log(\om^{-\ell}z), \frac{(\log(\om^{-\ell}z))^2}{2!}, \ldots,  (-1)^{a-1} \frac{(\log(\om^{-\ell}z))^{a-1}}{(a-1)!},0 , 0,\ldots,0\Big)$$
where the coefficient $-1$ in $Y_{0,\ell}(z) $ (resp. $1 $ in $Y_{\infty,\ell}(z) $) is in position $a+1+\ell$. 

We let $J_0 = \{(0,1),(0,2),\ldots,(0,N)\}$, $J_\infty = \{(\infty,1),(\infty,2),\ldots,(\infty,N)\}$, $J_{\om^\ell}  = \{\om^\ell\}$ for $1\leq \ell \leq N$, and  $\Sigma = \{0,\infty\} \cup \{\om^\ell, 1\leq \ell \leq N\}$, so that we have a solution $Y_j$ for each $j\in J_\sigma$, $\sigma\in\Sigma$.  

We also let $P_{a+1}(z) = \Pbar(z)$ (which is equal to $ \Pbar_{0,\ell}(z) + \Pbar_{\infty,\ell} (z)= (U_1+V_1)(\om^\ell z)$ for any $\ell$),  and $P_{a+1+\ell}(z) = \Pbar_{\infty,\ell}(z) = V_1(\om^\ell z) $ for any $\ell\in\unN$. Then 
we have polynomials $P_1(z)$, \ldots, $P_q(z)$ of degree at most $n$, and  with the notation of \S \ref{subsecshid} the remainders associated with the local solutions $Y_j$, $j\in J_\sigma$, $\sigma\in\Sigma$, are exactly the functions that appear in the Pad\'e approximation problem of Step 1:   $R(Y_{0,\ell}) = R_{0,\ell}(z)$, $R(Y_{\infty,\ell}) = R_{\infty,\ell}(z)$, and $R(Y_{\om^\ell}) = R_{\om^\ell}(z)$  for any $\ell\in\unN$.

Since $P_a$ is not the zero polynomial, we have $R_{\om^\ell}(z)\neq 0$ for any $\ell$; the functions  $R_{0,1}(z)$, \ldots, $R_{0,N}(z)$  (resp.  $R_{\infty,1}(z)$, \ldots, $R_{\infty,N}(z)$) are $\C$-linearly independent  (see  \cite[Lemma 2]{SFcaract}).
 
Eqns.  \eqref{eqP2},   \eqref{eqP3} and \eqref{eqP1bis} yield  $\ord_\infty (R_{\infty,\ell}(z)) \geq rn+1$,   
 $\ord_{\om^\ell}  (R_{\om^\ell}(z)) \geq \dz-1$ and  $\ord_0 (R_{0,\ell}(z)) \geq (r+1)n+1$  for any $\ell\in\unN$, so that 
 $$\sum_{\sigma\in\Sigma} \sum_{j\in J_\sigma} \ord_\sigma R_j (z) \geq (2r+1)Nn+N(\dz+1) = (n+1) q-nN -\tau \mbox{ with } \tau =a+1-aN; $$
 here $q = a+N+1$, and we recall that $\dz = -\deg F = (a+1)  (\frac{n}{N} +1) - (2r+1)n -1$. As above,  $\tau$ is exactly the 
 difference between the number of unknowns and the number of equations
in the Pad\'e approximation problem of Step 1.  

The definition \eqref{eqdefpkishid} of $P_{k,i}$ is consistent with the one given (for $i\leq a$) by Eq. \eqref{eqdefpkj} in \S \ref{subsec33}. We have $\tau = a+1-aN$, so that for $n$ sufficiently large $\tau+ \cdeux \leq \dz$ where $\cdeux$ is the constant given by Theorem \ref{thshid}. Therefore \eqref{eqokenun} shows that $U_k$ and $V_k$ are holomorphic at $z=1$ for any $k < \tau+ \cdeux$. Eqns. \eqref{eqdefpkj},  \eqref{eqdefpz} and  \eqref{eqdefpinf} imply that they are holomorphic at all other roots of unity. Now  Eqns. \eqref{eqdefpkishid},  \eqref{eqdefpz} and  \eqref{eqdefpinf} yield
\begin{equation}\label{eqnvp}
P_{k,a+1}(z)=\om^{\ell(k-1)} ( U_k(\om^\ell z) + V_k ( \om^{\ell} z)) \mbox{ and } 
P_{k,a+1+\ell}(z) = \om^{\ell(k-1)}  V_k ( \om^{\ell} z) 
\eneq
for any $\ell \in\{1,\ldots,N\}$, by induction on  $k\geq 1$. Therefore all $P_{k,i}$, with $k < \tau+\cdeux$ and $1\leq i \leq q$, are holomorphic at 1. 

\bigskip

We have checked all assumptions of Theorem \ref{thshid} for $n$ sufficiently large: the matrix $[P_{k,i}(1)]_{1\leq i \leq q, 1\leq k < \tau+\cdeux}$ has rank at least $q-\Card J_1 =q-1$. Now \eqref{eqokenun} implies $P_{k,1}(1) = 0$ for any $k < \tau+\cdeux$, so that we may remove the first row: the matrix $P = [P_{k,i}(1)]_{2\leq i \leq q, 1\leq k < \tau+\cdeux}$ has rank $q-1$, equal to its number of rows. 
\bigskip

\noindent {\bf Step 3:} Expression of $s_{k,i}$ in terms of $P$ and conclusion.

\bigskip

We shall now compute a matrix $M$ {\em independent from $n$} such that $[s_{k,i}]_{i,k} = MP$; recall that the coefficients $s_{k,i}$ and the matrix $P$ depend on $n$.

To begin with, Eq. \eqref{eqpkzinfom} with $z = \om^\ell$ yields
$$
\om^{ (k-1)\ell } U_k(\om^\ell ) = \sum_{\lam = 0}^{N-1} \om^{\lam\ell}  U_{k,\lam}(1) \mbox{ for any } \ell\in\Z,$$
since $U_{k,\lam}(z) \in \Q[z^N, z^{-N}]$. Therefore we have 
\begin{equation}\label{eqdirecteun}
  U_{k,\lam}(1) = \frac1N \sum_{\ell=1}^N \om^{(k-1-\lam)\ell} U_k(\om^\ell) \mbox{ for any } 0\leq\lam\leq N-1,
  \eneq
and the same relation holds with $V_{k,\lam}$ and $V_k$. Using Eq. \eqref{eqnvp} we deduce that 
$$V_{k,\lam}(1) =  \frac1N \sum_{\ell=1}^N \om^{-\lam\ell} P_{k, a+1+\ell} (1) \mbox{ for } 0\leq\lam\leq N-1,$$
and also for $\lam=N$ since $V_{k,N} = V_{k,0}$, and 
$$U_{k,\lam}(1) =
\left\{\begin{array}{l}
 - \frac1N \sum_{\ell=1}^N \om^{-\lam\ell} P_{k, a+1+\ell} (1) \mbox{ if } 1\leq\lam\leq N-1, \\
P_{k, a+1} (1) - \frac1N \sum_{\ell=1}^N  P_{k, a+1+\ell} (1) \mbox{ if } \lam=0. \\
\end{array}\right.
$$
Therefore the definition  \eqref{eqdefski} of $s_{k,i}$ can be translated as 
\begin{equation}\label{eqpromatexpli}
s_{k,i}  = \sum_{j=2}^q m_{i,j} P_{k,j}(1)
\eneq
for any $2\leq i \leq a+N$ and any $1\leq k \leq \dz-1$, where  the coefficients $m_{i,j}$ are defined for $2\leq i \leq a+N$ and $2\leq j \leq q=a+N+1$ by 
\begin{equation}\label{eqdefm}
\left\{\begin{array}{l}
m_{i,i} = \delta_n \mbox{ for  } 2\leq i \leq a+1 \\
m_{a+1, a+1+\ell} = \frac{\delta_n}{N} ((-1)^p-1)   \mbox{ for  } 1\leq \ell \leq N \\
m_{a+1+\lam, a+1+\ell} = \frac{\delta_n}{N} ((-1)^p \om^{\lam \ell} - \om^{-\lam \ell})   \mbox{ for  } 1\leq \ell \leq N \mbox{ and } 1\leq \lam \leq N-1 \\
m_{i,j} = 0 \mbox{ for all other pairs } (i,j).
\end{array}\right.
\eneq
Let us choose now the constant $\cstnum$  of Lemma \ref{lemmatinv}; the same constant appears in Proposition \ref{propdio}. We take $ \cstnum = \tau+\cdeux -1 $; this constant depends only on $a$ and $N$. We consider the matrices  $M = [ m_{i,j}  ]_{2\leq i \leq a+N, 2\leq j \leq q}$ and $P = [P_{k,j}(1)]_{2\leq j \leq q, 1\leq k \leq \cstnum}$. Then  Eq. \eqref{eqpromatexpli} means that 
\begin{equation}\label{eqpromat}
[s_{k,i}]_{2\leq i \leq a+N, 1\leq k \leq \cstnum } = MP.
\eneq
Both $M$ and $P$  have coefficients in $\Q(\om)$; recall that the coefficients $s_{k,i}$ of $MP$ are rational numbers, and we shall prove in \S \ref{subsec35} that they are integers. Let $\calF$ denote the subspace of $\R^{a+N-1}$ spanned by the $q-1$ columns $\tra ( m_{2,j}, \ldots, m_{a+N,j})$ of $M$. Now assume that $n$ is sufficiently large; then we have proved in Step 2 that  the $q-1$ rows of $P$ are linearly independent. Therefore  Eq. \eqref{eqpromat} shows that $\calF$ is equal to the subspace spanned by columns $ \tra (s_{k,2},\ldots, s_{k, a+N})$ of the matrix $[s_{k,i}]_{i,k} $. Since $M$ does not depend on $n$, neither does $\calF$: this concludes the proof of  Lemma \ref{lemmatinv}.

\bigskip

\begin{Rem} \label{remf} Let us prove that in Lemma \ref{lemmatinv}, the subspace $\calF$ is not always equal to  $\R^{a+N-1}$ (i.e., that the matrix $[s_{k,i}]$ may have rank less than its number of rows, namely $a+N-1$). Consider the case where $p$ and $N$ are even (so that $\om^{N/2} = -1$). Then  the definition \eqref{eqdefm} of the matrix $M$ in Step 3 above yields $m_{a+1+N/2, j} = 0$ for any $j$, so that Eq. \eqref{eqpromatexpli} implies $s_{k, a+1+N/2}=0$ for any $k$: the matrix $[s_{k,i}]_{i,k}$ has a zero row. This phenomenon does not occur in \cite{SFcaract}; it comes from the new property \eqref{equplusv} obtained in Remark \ref{remnv}. Indeed a direct proof that $s_{k, a+1+N/2}=0$ can be obtained as follows, 
 using Eqns. \eqref{eqdefski}, \eqref{eqdirecteun}, and \eqref{eqnvp} but not 
  the matrix $M$: 
\begin{eqnarray*}
s_{k, a+1+N/2}&=&\delta_n \Big( U_{k,N/2}(1) + V_{k,N/2}(1) \Big) \\
&=& \frac{\delta_n }{N} \sum_{\ell=1}^N \om^{(k-1-N/2)\ell}(U_k+V_k)(\om^\ell) = 
 \frac{\delta_n }{N} P_{k,a+1}(1) \sum_{\ell=1}^N (-1)^\ell = 0.
 \end{eqnarray*}
\end{Rem}

\subsection{Arithmetic and Asymptotic Properties}\label{subsec35}

In this section we conclude the proof of Proposition~\ref{propdio} stated in  \S \ref{subsec21}, by proving parts $(ii)$ and $(iii)$ and the fact that the $s_{k,i}$ are integers. Recall that 
$$
\alpha = (4e)^{(a+1)/N} (2N)^{2r+2} r^{-(a+1)/N+4(r+1)} \mbox{ and }
\beta = (2e )^{(a+1)/N} (r +1)^{2r+2} N^{2r+2}.
$$

\begin{Lem} \label{lemtech}
We have $s_{k,i}\in\Z$ for any $i\in\{2,\ldots,a+N\}$ and any $k\leq \dz-1$, and as $n\to\infty$:
$$\Big| 2(-1)^p \sum_{\substack{i=2 \\ i\equiv p \bmod 2}}^a s_{k,i} \xii + \sum_{i=0}^{N-1} s_{k,a+1+i}f(i)\Big| \leq \alpha ^{n+o(n)}
\mbox{ and }
 \max_{2\leq i \leq a+N}  | s_{k,i} | \leq \beta  ^{n(1+o(1))}.$$
\end{Lem}

In this lemma and throughout this section, we denote by $o(1)$ any sequence that tends to 0 as $n\to\infty$; it usually depends  also on $a$, $r$, $N$, and $k$ (but the dependence in $k$ is not significant since $k$ is bounded by $d_0-1$, which depends only on $n$, $a$, $r$, $N$). We also   recall  that $d_{n}$ is the least common multiple of 1, 2, \ldots, $n$.

\bigskip

We shall prove two lemmas now; the deduction of Lemma \ref{lemtech} from these lemmas (using  Lemma \ref{lemdcp} proved in \S \ref{subsec33})  is exactly the same as the proof of Proposition 1 in \cite[\S 4.5]{SFcaract}. Recall  from \S \ref{subsec33} that 
$$F(t) = (n/N)!^{(a+1)-(2r+1)N} \frac{(t-rn)_{(2r+1)n+1}  }{\prod_{h=0}^{n/N}(t+Nh)^{a+1}} 
  = \sum_{h=0}^{n/N}\sum_{j=1}^a \frac{p_{j,h}}{(t+Nh)^j}. $$

\begin{Lem} \label{lempjh}
For any $j\in\una$ and any $h \in \zeronsN$ we have
\begin{equation} \label{eqdenompjh}
(N d_{n/N})^{a+1-j} N^{(a+1)n/N} p_{j,h} \in \Z
\eneq
\begin{equation} \label{eqmajopjh}
\mbox{ and }
| p_{j,h} | \leq \Big(   2^{(a+1)/N}N^{2(r+1)-(a+1)/N}(r+1)^{2r+2}     \Big)^{n(1+o(1))}
\eneq
where $o(1)$ is a sequence that tends to 0 as $n\to\infty$ and may depend also on $N$, $a$, and $r$.
\end{Lem}

\Dem of Lemma \ref{lempjh}:   We follow the approach of   \cite{Colmez} by letting
\begin{eqnarray*}
F_0(t)  &= \frac{(n/N)!}{\prod_{h=0}^{n/N} (t+Nh) } & = \sum_{h=0}^{n/N} \frac{(-1)^h N^{-n/N} \combitiny{n/N}{h}}{t+Nh} ,\\ 
G_i(t)  &= \frac{(t-in/N)_{n/N}}{\prod_{h=0}^{n/N} (t+Nh) } & = \sum_{h=0}^{n/N} \frac{(-1)^{h+n/N} N^{-n/N} \combitiny{n/N}{h}  \combitiny{Nh+in/N}{n/N}}{t+Nh}  \mbox{ for } 1\leq i \leq rN, \\ 
H_i(t)  &= \frac{(t+1+in/N)_{n/N}}{\prod_{h=0}^{n/N} (t+Nh) } & = \sum_{h=0}^{n/N} \frac{(-1)^{h } N^{-n/N} \combitiny{n/N}{h}  \combitiny{-Nh+(i+1)n/N}{n/N}}{t+Nh}  \mbox{ for } 0\leq i \leq (r+1)N-1. 
\end{eqnarray*}
Then the partial fraction expansion of 
$$F(t)= F_0(t) ^{a+1-(2r+1)N} t G_1(t)\ldots G_{rN}(t) H_0(t) \ldots H_{(r+1)N-1}(t)$$
 can be obtained by multiplying those of $F_0$, $G_i$ and $H_i$ using repeatedly the formulas $\frac{t}{t+Nh} = 1 - \frac{Nh}{t+Nh}$ and 
\begin{equation} \label{eqprodelsples}
\frac1{(t+Nh)(t+Nh')^\ell } = \frac1{N^\ell (h'-h)^\ell (t+Nh)} - \sum_{i=1}^\ell \frac1{N^{\ell+1-i} (h'-h)^{\ell+1-i} (t+Nh')^i} 
\eneq
with $h\neq h'$. The denominator of $p_{j,h} $ comes both from this formula (and this contribution divides $(Nd_{n/N})^{a+1-j}$) and from the denominators of the coefficients in the partial fraction expansions of $F_0$, $G_i$, $H_i$ (which belong to $N^{-n/N}\Z$, so that $N^{(a+1)n/N}$ accounts for this contribution). This concludes the proof of \eqref{eqdenompjh}.

\bigskip

On the other hand, bounding from above the coefficients of the 
partial fraction expansions of $F_0$, $G_i$, $H_i$  yields
$$|p_{j,h}| \leq n^{O(1)} N^{-(a+1)n/N} 2^{(a+1)n/N}\prod_{i=1}^{rN } \frac{(n+in/N)!}{(n/N)! (n+(i-1)n/N)!} \prod_{i=0}^{(r+1)N-1 } \frac{((i+1)n/N)!}{(n/N)! ( i n/N)!} $$
where $O(1)$ is a constant depending only on $a$, $r$, $N$ which can be made explicit (see \cite{Colmez} for details). Simplifying the products and using the bound $\frac{m!}{m_1!\ldots m_c!} \leq c^m$ valid when $m_1+\ldots+m_c = m$, one obtains
$$|p_{j,h}| \leq n^{O(1)} (2/N)^{(a+1)n/N}  \frac{((r+1)n)!^2}{n! (n/N)!^{(2r+1)N}}  \leq  n^{O(1)} (2/N)^{(a+1)n/N} ((r +1)N) ^{2(r+1)n}.$$
This concludes the proof of Lemma \ref{lempjh}.

\bigskip

The proof of the  following lemma is inspired by that of \cite[Lemma 1.4]{Sprang}. Recall that $U_1$ and $V_1$ are defined in Eqns. \eqref{eqdefu1} and \eqref{eqdefv1}, and that 
$$\delta_n = (N d_{n/N})^{a+1} N^{(a+1)n/N}.$$

\begin{Lem} \label{lemsprang} The polynomials $\delta_n U_1(z)$ and  $\delta_n V_1(z)$ have integer coefficients.
\end{Lem}

\Dem of Lemma \ref{lemsprang}: Recall from Eq. \eqref{eqdefu1} that 
$$U_1(z) =  -\sum_{t=1}^n z^t \sum_{j=1}^a \sum_{h=0}^{\lfloor (t-1)/N\rfloor} \frac{p_{j,h}}{(Nh-t)^j}  . $$
Assume that  $\delta_n U_1(z)$ does not have integer coefficients. Then there exists $t\in\{1,\ldots,n\}$ such that 
$$\sigma :=  \sum_{j=1}^a \sum_{h=0}^{\lfloor (t-1)/N\rfloor} \frac{ (N d_{n/N})^{a+1} N^{(a+1)n/N} p_{j,h}}{(Nh-t)^j} \not\in\Z.$$
Let $p'_{j,h} = (N d_{n/N})^{a+1-j} N^{(a+1)n/N} p_{j,h}$, which is an integer thanks to   Lemma \ref{lempjh}. Then we have
$$\sigma =  \sum_{j=1}^a \sum_{h=0}^{\lfloor (t-1)/N\rfloor} \frac{   d_{n/N}^j p'_{j,h} }{ (h-t/N)^j}.$$
If $N$ divides $t$ then $\frac{t}{N}-h$ is a positive integer less than or equal to $n/N$, so that it divides $d_{n/N}$: this contradicts the assumption $\sigma\not\in\Z$. Therefore $N$ does not divide $t$, so that $F(-t)=0$.

Since  $\sigma\not\in\Z$ there exists $h_0$ such that 
$$ \sum_{j=1}^a \frac{   d_{n/N}^j p'_{j,h_0} }{(h_0-t/N)^j} \not\in\Z.$$
Now $F(-t)=0$ so that 
$$  \sum_{j=1}^a \sum_{\substack{ h=0 \\ h\neq h_0}}^{n/N}  \frac{   d_{n/N}^j p'_{j,h} }{ (h-t/N)^j}
= -  \sum_{j=1}^a \frac{   d_{n/N}^j p'_{j,h_0} }{(h_0-t/N)^j} \not\in\Z.$$
This rational number has negative $p$-adic valuation for some prime number $p$. Therefore there exist $h_1\neq h_0$ and $j_0$, $j_1$ such that 
$$v_p \Big( \frac{d_{n/N}^{j_1}}{(h_1 - t/N)^{j_1}}\Big)< 0 \mbox{ and } 
v_p \Big( \frac{d_{n/N}^{j_0}}{(h_0 - t/N)^{j_0}}\Big)< 0.$$
This implies
$$\min ( v_p ( h_1-t/N),  v_p ( h_0-t/N)) > v_p ( d_{n/N})$$
so that $v_p (h_0 - h_1) > v_p ( d_{n/N})$. This is a contradiction since $1\leq | h_0 - h_1| \leq n/N$. This concludes the proof that $\delta_n U_1(z) \in \Z[z]$; the same proof works for  $\delta_n V_1(z) $.

\section{Siegel's linear independence criterion}\label{sec3}

The following criterion  is based on  Siegel's ideas (see for instance \cite[p. 81--82 and 215--216]{EMS}, \cite[\S 3]{Matala-Aho}  or  \cite[Proposition 4.1]{Marcovecchio}).

\begin{Prop} \label{propsiegel}
Let $\theta_1,\ldots,\theta_q$ be  complex numbers, not all zero. Let $\tau>0$, and $(Q_n)$ be a sequence of real numbers with limit $+\infty$. Let $\calN$ be an infinite subset of $\N$, $K\geq 1$, and for any $n\in\calN$ let $L^{(n)} = [\ell_{k,i}^{(n)}]_{1\leq i\leq q, 1\leq k \leq K}$ be a matrix with integer coefficients  such that as $n\to\infty$ with $n\in\calN$:
$$\max_{i,k  } | \ell_{k,i}^{(n)}| \leq Q_n^{1+o(1)}$$
$$
\mbox{ and } 
\max_{1\leq k \leq K } | \ell_{k,1}^{(n)} \theta_1 + \ldots  + \ell_{k,q}^{(n)} \theta_q | \leq Q_n^{-\tau + o(1)}.$$
Assume also that the subspace $\calF$ of $\R^q$ spanned by the columns $\tra (\ell_{k,1}^{(n)}, \ldots, \ell_{k,q}^{(n)} )$ of $L^{(n)}$ is non-zero and independent from $n\in\calN$ (provided $n$ is large enough).  Then we have
$$\dim_\Q\Span_\Q(\theta_1,\ldots,\theta_q)\geq \tau+1.$$
\end{Prop}

The usual version of this criterion (see for instance  \cite[Theorem 4]{gfndio2}) is the same statement, but the assumption on $\calF$ is replaced by the assumption that    $L^{(n)}$  is   invertible. The latter is stronger, since it  is equivalent to asking $\calF = \R^q$ for any $n$.  Indeed if $\calF=\R^q$ then $L^{(n)}$ has rank $q$: for each $n$ we may extract $q$ linearly independent columns of $L^{(n)}$, and obtain an invertible matrix to which    \cite[Theorem 4]{gfndio2} applies. The point is that we shall apply Proposition \ref{propsiegel} to the matrices $[s_{k,i}]$ constructed in Proposition \ref{propdio}, and the subspace $\calF$ is not always equal to $\R^q$ (see Remark  \ref{remf} in \S \ref{subsec34}). 

\bigskip

Let us prove  Proposition \ref{propsiegel} now. Denote by $\calF$ the image of  $L^{(n)}$, assumed to be independent from $n \in\calN$ (provided $n$ is large enough). Let $p=\dim \calF$. Permuting $\theta_1$, \ldots, $\theta_q$ if necessary, we may assume that a system of linear equations of $\calF$ is given by
\begin{equation} \label{eqff}
x_t = \sum_{i=1}^p \mu_{t,i} x_i \mbox{ for } p+1\leq t \leq q, \mbox{ with } \mu_{t,i}\in\Q.
\eneq
We point out that the coefficients $\mu_{t,i}$ are rational numbers because the matrices $L^{(n)}$ have integer coefficients. Since $\tra ( \ell_{k,1}^{(n)}, \ldots, \ell_{k,q}^{(n)} )\in\calF $ for any $1\leq k \leq K$ and any $n \in\calN$ sufficiently large, Eq. \eqref{eqff} yields
\begin{equation} \label{eqtheta}
\sum_{i=1}^p  \ell_{k,i} ^{(n)}   \theta_i = \sum_{i=1}^p  \ell_{k,i} ^{(n)}  \Big(  \theta_i  +  \sum_{t=p+1}^q \mu_{t,i}  \theta_t\Big)
=  \sum_{i=1}^p  \ell_{k,i} ^{(n)}    \theta'_i  
\eneq
upon letting $ \theta'_i  =   \theta_i  +  \sum_{t=p+1}^q \mu_{t,i}  \theta_t$ for $1\leq i \leq p$. Moreover for any $n \in\calN$ sufficiently large, we have $\rk L^{(n)} = \dim \calF = p$ and Eq. \eqref{eqff} shows that the last $q-p$ rows of $L^{(n)}$  are linear combinations of the first $p$ rows. Therefore the first $p$ rows are linearly independent: the matrix $[ \ell_{k,i}^{(n)}]_{1\leq i\leq p, 1\leq k \leq K}$ has rank $p$. Accordingly for each $n$ there exist pairwise distinct integers $k_1$, \ldots, $k_p$ between 1 and $K$ such that the matrix $M^{(n)} = [ \ell_{k_j,i}^{(n)}]_{1\leq i,j\leq p}$   is invertible. Using Eq. \eqref{eqtheta} we may apply the usual version of Siegel's criterion (namely    \cite[Theorem~4]{gfndio2}) to this matrix and deduce that 
$$\dim_\Q\Span_\Q(\theta'_1,\ldots,\theta'_p)\geq \tau+1.$$
Since $\theta'_i \in \Span_\Q(\theta_1,\ldots,\theta_q)$ for any $1\leq i \leq p$, this concludes the proof of Proposition \ref{propsiegel}.

\bigskip

\begin{Rem} \label{remgfndio} The idea of applying the usual version of Siegel's criterion to numbers $\theta'_i$ defined as linear combinations of $\theta_1$, \ldots, $\theta_q$ appears also in \cite{gfndio2} (see Proposition 2 in \S 6 and Eq. (9.1)). However the situation is different in that paper: the rows of the matrix $P$ (see Step 2 in \S \ref{subsec34} above) are linearly dependent, which is not the case here.
\end{Rem}

\section{Deduction of Theorems \ref{th1}  and \ref{th2} from Proposition \ref{propdio}}\label{sec1}

In this section we prove Theorems \ref{th1}  and \ref{th2} stated in the introduction, and also a result that nearly contains both of them (namely Theorem \ref{th3} stated at the end of \S \ref{subsecdemth2}).  At last, we show in \S \ref{subsecFSZ} that the linear forms constructed in \cite{FSZ} are a special case of those studied here.

\subsection{Proof of Theorem  \ref{th1}} \label{subsecdemth1}

Let $f$,  $T$, $p$, $\eps$, and $a$  be as in the statement of Theorem  \ref{th1}; put $N=T$. We consider the complex numbers $\theta_1,\ldots, \theta_{a+N-1}$ given by:
$$
\left\{ 
\begin{array}{l}
\theta_{i-1} = 2(-1)^p L(f,i) \mbox{ for $2\leq i \leq a$ with $i\equiv p \bmod 2$,}\\
\theta_{i-1} = 0 \mbox{ for $2\leq i \leq a$ with $i\not\equiv p \bmod 2$,}\\
\theta_{a+i} = f(i)  \mbox{ for } 0\leq i \leq N-1 .
\end{array}\right. $$
We apply Proposition \ref{propdio} to each integer multiple $n$ of $N$, and let $ \ell_{k,i}^{(n)}   = s_{k,i+1}$ for $1\leq i \leq a+N-1$ and $1\leq k \leq \cstnum$. Then we   apply Siegel's linear independence criterion  (namely Proposition \ref{propsiegel} stated and proved in  \S \ref{sec3}) with $q=a+N-1$, $Q_n = \beta^n$ and $\tau = -\frac{\log \alpha}{\log\beta}$ (so that $Q_n^{-\tau} = \alpha^n$), where $\alpha $ and $\beta$ are defined in \S \ref{subsec21};  we take for $\calN$  the set of integer multiples of $N$. Therefore we obtain 
\begin{equation} \label{eqminopreuveun}
\dim_\Q\Span_\Q\Big(    \{ \xii, \, 2\leq i \leq a \mbox{ and  } i\equiv p \bmod 2 \}
\cup   \{ f(0) , \ldots,  f(N-1)\} 
 \Big)\geq 1 - \frac{\log\alpha}{\log\beta}.
 \end{equation}
 Taking  $a$ very large  and $r$ equal to the integer part of $\frac{a}{(\log(a))^2}$  concludes  the proof of Theorem  \ref{th1} since 
$$1 - \frac{\log\alpha}{\log\beta} -N = \frac{1+\eps_a}{1+\log 2} \log a \mbox{ where } \lim_{a\to+\infty} \eps_a = 0;$$
here the  shift of $N$ in the lower bound comes from $f(0)$, \ldots, $f(N-1)$ that  appear in Eq.~\eqref{eqminopreuveun}.

\subsection{Proof of Theorem  \ref{th2}} \label{subsecdemth2}

Let $f$,  $T$, $E$, and  $p$  be as in the statement of Theorem  \ref{th2}. Let $0 < \eps < 1/4$ and $a $ be  sufficiently large with respect to $\eps$, $T$, and $\dim E$. We denote by $D$   the product  of all primes less than or equal to $(1-3\eps) \log a$ (such a product has asymptotically the largest possible number of divisors with respect to its size, see \cite[Chapter XVIII, \S 1]{HW}). Then we have
$$\log D \, \,  = \sum_{p\leq  (1-3\eps) \log a} \log p \leq ( 1-2\eps) \log a$$
by the prime number theorem, i.e.,  $D\leq   a^{1-2\eps}$. We take for $r$ the integer part of $a^{\eps}$. At last, we let $N = DT$.

\bigskip

For any divisor $d$ of $D = N/T$ and any $m\in\Z$, let $g_d(m ) = f(m/d)$ if $d$ divides $m$, and $g_d(m ) =0$ otherwise. Since $f$ is $T$-periodic we have $g_d(m+N) = g_d(m)$ for any $m$.

We shall  choose below an   integer $w_d$ for each divisor $d$ of $D$; let 
$$g = \sum_{d | D} w_d g_{D/d}.$$

We shall apply  Proposition \ref{propdio}  to the $N$-periodic function $g$ and obtain   linear forms in the numbers 
$$L(g,i)= \sum_{m=1}^{\infty} \frac1{m^i} \sum_{\substack{d | D \\  D|md }} w_d f(md/D) = 
\sum_{ d | D } w_d \sum_{m'\geq 1} \frac{f(m')d^i}{ m'^i D^i}$$
by letting $m'= md/D$. Therefore we have
\begin{equation} \label{eqxii}
L(g,i)= D^{-i}  \Big( \sum_{ d | D} w_d d^{ i}\Big)L(f,i).
 \end{equation}

\bigskip

Notice that $D$ has $\delta = 2^{\pi( (1-3\eps) \log a)}$ divisors, with 
\begin{equation} \label{eqamodif}
\log \delta = \pi( (1-3\eps) \log a) \log 2 \geq (1-4\eps) (\log 2) \frac{\log a}{\log \log a}.
\eneq
Assume that the number of values \eqref{eq01} which do not belong to $E$ is   less than $\delta$. Let $2 \leq  i_1 < i_2 < \ldots < i_{\delta-1}\leq a$ be  integers such that if $ L(f,i) \not\in E$ and $i\equiv p \bmod 2$, $2\leq i \leq a$, 
  then $i=i_j$ for some $j$.
  
The homogeneous linear system
\begin{equation} \label{eqwun}
\sum_{d | D}w_d d^{ i_j} = 0 \mbox{ for any } j\in\{1,\ldots,\delta-1\}
\end{equation}
has $\delta$ unknowns $w_d$, where $d$ ranges through the set $\calD$ of divisors of $D$, and $\delta-1$ equations. Therefore it has a non-zero integer solution $(w_d)\in\Z^\calD$. 

At this point, the integers $w_d$ are chosen in \cite{FSZ} such that $\sum_{d| D} w_d d \neq 0$, using an elementary zero estimate (namely, a generalized Vandermonde determinant is non-zero). Here we do not need to make any such assumption: we just assume that $w_d\neq 0$ for at least one $d$. Indeed a (much more complicated) zero estimate is used in the present proof, namely Theorem \ref{thshid}.

 Proposition \ref{propdio} applies to the $N$-periodic function $g = \sum_{d | D} w_d g_{D/d} $ defined above. Using also Siegel's linear independence criterion as in \S \ref{subsecdemth1} we obtain
\begin{equation} \label{eqdimmino}
\dim_\Q\Span_\Q\Big( \{ g(0) , \ldots,  g(N-1)\} \cup \{ L(g,i), \, 2\leq i \leq a \mbox{ and  } i\equiv p \bmod 2 \}\Big)\geq  1 - \frac{\log\alpha}{\log\beta}  
\end{equation}
with
 $$1 - \frac{\log\alpha}{\log\beta}  \sim\frac{\log r}{1+\log 2} \sim \frac{\eps}{1+\log 2}\log a $$
as $a\to\infty$ (recall that $r$ is  the integer part of $a^{\eps}$). 

On the other hand, the numbers that appear in the left hand side of \eqref{eqdimmino} have the following properties:
\begin{itemize}
\item[$\bullet$] $g(0)$, \ldots, $g(N-1)$  belong to $\{0,f(0),f(1),\ldots,  f(T-1)\}$.
\item[$\bullet$] For $2\leq i \leq a$ with   $i\equiv p \bmod 2$, $L(g,i)$  is zero if $i\in\{i_1,\ldots,i_{\delta-1}\}$, and belongs to $E$ otherwise, 
 as Eqns. \eqref{eqxii} and \eqref{eqwun} show.
 \end{itemize}
Therefore we have 
\begin{equation} \label{eqdimmajo}
\dim_\Q\Span_\Q\Big( \{ g(0) , \ldots,  g(N-1)\} \cup \{ L(g,i), \, 2\leq i \leq a \mbox{ and  } i\equiv p \bmod 2 \}\Big)\leq T +\dim E.
\end{equation}
Combining Eqns. \eqref{eqdimmino} and  \eqref{eqdimmajo} yields a contradiction   provided $a$ is large enough. This   concludes the proof of Theorem  \ref{th2}.

\bigskip

Since $4(1+\log 2) > 7$, the same proof (with $\eps$ replaced with $\eps/4$ to take Eq. \eqref{eqamodif} into account) provides the following refinement of Theorem \ref{th2}.
 
 \begin{Th} \label{th3}
Let $T \geq 1$, and $f: \Z \to \C$ be such that $f(n+T)=f(n)$ for any $n$. Assume that $f$ is not identically zero. Let $p\in\{0,1\}$, $0< \eps < 1$, and $a $ be sufficiently large (in terms of  $T$  and $\eps$).   Let $E$ be a finite-dimensional $\Q$-vector space contained in $\C$ with $\dim E < \frac{\eps}{7} \log a$.  Then among the numbers $L(f,s)$   with $2\leq s \leq a$ and $s\equiv p \bmod 2$,
at least
$$2^{(1-\eps)\frac{\log a}{\log \log a}}$$
do not belong to $E$.
\end{Th}

Choosing $\eps = 7/8$,  this refinement implies that the numbers $L(f,s)$   with $2\leq s \leq a$ and $s\equiv p \bmod 2$ are not all contained in such a subspace $E$: they span a $\Q$-vector space of dimension at least $\frac18 \log a$. Except for the multiplicative constant ($\frac18$ instead of $\frac{1-o(1)}{1+\log 2}$),  Theorem \ref{th1} follows as a corollary of Theorem \ref{th3}.

\subsection{Connection to the proof of \cite{FSZ}}\label{subsecFSZ}

In this section we show that the linear forms used in \cite{FSZ} to prove \eqref{eqFSZ} are a special case of those studied in the present paper (namely in the proof of Theorem \ref{th2} with $f(m)=1$, $T=r=k=1$, and $p\equiv a \bmod 2$). Accordingly they are related to the  Pad\'e approximation problem stated in \S \ref{subsec34}, in which the number of equations is essentially equal to the number of unknowns.

\bigskip

We keep the notation of the proof of  Theorem \ref{th2} in \S \ref{subsecdemth2}, with $T=1$ and $f(m)=m$ for any $m\in\Z$. In particular $N=D$ is   the product  of all primes less than or equal to $(1-3\eps) \log a$. For any divisor $d$ of $D$ we have $g_d(m)=1$ if $d$ divides $m$, and $g_d(m)=0$ otherwise. The function $g = \sum_{d|D} w_d g_{D/d}$ satisfies
\begin{equation} \label{eqgg}
g(m) =    \sum_{\substack{d | D \\  D|md }} w_d \mbox{ for any } m\in\Z.
\eneq
Now let $n$ be an integer multiple of $2N=2D$, and let $p\in\{0,1\}$ be such that $p\equiv a \bmod 2$.  Then the rational function $F$ satisfies the symmetry property of well-poised hypergeometric series:
\begin{equation} \label{eqwp}
F(-n-t) = (-1)^{(2r+1)n+1+(a+1)(\frac{n}{N}+1)}F(t) = (-1)^p F(t).
\eneq
This is the key ingredient (since the Ball-Rivoal theorem) to get rid of even zeta values, when $p=1$. In our approach where Nesterenko's linear independence criterion is replaced with Siegel's combined with Shidlovsky's lemma, this property cannot be used in the same way because it is destroyed when considering $S^{(k-1)}(z)$ for $k\geq 2$. Using both $S_0$ and $S_\infty$ in constructing the linear forms (see \S \ref{subsec33}) makes it possible to overcome this difficulty (as in \cite{SFcaract}). With $k=1$ this trick does not modify the linear forms we are interested in, since for any $\ell\in\Z$ we have using Eqns. \eqref{eqdefszinf} and  \eqref{eqwp} and the fact that $N$ divides $n$:
\begin{eqnarray}
S_0(\om^\ell)+(-1)^p S_\infty(\om^{-\ell})
&=& 
\sum_{t=n+1}^{\infty} F(-t) \om^{ \ell t} + (-1)^p \sum_{t=1}^{\infty} F( t) \om^{ \ell t}  \nonumber   \\ 
&=& 
\sum_{t=1}^{\infty}   \om^{ \ell t}  \Big( F(-n-t) + (-1)^p  F( t) \Big)    \nonumber   \\ 
 &=& 
2(-1)^p \sum_{t=1} ^{\infty} F( t)  \om^{ \ell t} . \label{eqsimpli}
\end{eqnarray}

\bigskip

We are now in position to express differently the linear forms constructed in part $(iii)$ of Proposition \ref{propdio} from the map $g$   given by Eq. \eqref{eqgg}, in the special case where  $N=D$, $n$ is a multiple of $2N$, $p\equiv a \bmod 2$, $r=1$, and $k=1$. Denote by $\Lambda_n$ this linear form. Then we have using Lemma \ref{lemdcp} and Eqns. \eqref{eqsimpli} and \eqref{eqgg}:
\begin{eqnarray*}
\delta_n^{-1} \Lambda_n 
&=& \sum_{\ell=1}^D \widehat{g}(\ell)  \Big( S_0(\om^\ell)+(-1)^p S_\infty(\om^{-\ell}) \Big) \\
 &=& 
2(-1)^p \sum_{t=1} ^{\infty} F( t) \sum_{\ell = 1}^D  \widehat{g}(\ell) \om^{ \ell t} \\
 &=& 
2(-1)^p \sum_{t=1} ^{\infty} F( t) g( t) \\
 &=& 
2(-1)^p   \sum_{ d | D  } w_d   \sum_{\substack{t\geq 1 \\   D| d t  }}  F( t)   \\
 &=& 
2(-1)^p   \sum_{ d | D  } w_d   \sum_{t'=1} ^\infty   F( Dt'  / d)   \\
 &=& 
2(-1)^p   \sum_{ d | D  } w_d  \sum_{j=1} ^d  \sum_{m=1} ^\infty   F( mD + j \frac{D}{d}) .
\end{eqnarray*}
In the last expression the sum on $m\geq 1$ should have begun at $m=0$, but this makes no difference since $F(jD/d)=0$ for any $1\leq j \leq d$. Now let $R(t) = F(Dt)$; then we have 
\begin{equation} \label{eqlambdaeg}
\frac{ (-1)^p}{ 2\delta_n} \Lambda_n = 
  \sum_{ d | D  } w_d  \sum_{j=1} ^d  \sum_{m=1} ^\infty   R( m  +   \frac{j}{d}) . 
 \eneq  
Up  to the normalizing factor $\frac{ (-1)^p}{ 2\delta_n} $ these are exactly the linear forms $\widehat{ r_n}$ used in \cite{FSZ} to prove \eqref{eqFSZ}. Indeed the following notation is used in \cite{FSZ} for $1 \leq  j \leq D$ and $d| D$:
$$R_n(t) = D^{3Dn} \, \, n!^{s+1-3D} \, \, \frac{ \prod_{j=0}^{3Dn} (t-n+\frac{j}{D})}{ \prod_{j=0}^{ n} (t+j)^{a+1}}, \quad \quad r_{n,j} = \sum_{m=1}^\infty R_n\Big(m+\frac{j}{D}\Big),$$
$$\widehat r_{n,d} = \sum_{j=1}^d r_{n, j\frac{D}{d}},  \quad \quad
\widetilde r_n = \sum_{d | D} w_d \, \widehat r_{n,d}.$$
Now, up to the normalizing factor $R_n(t)$ is equal to the rational function $R(t) = F(Dt)$ so that $\Lambda_n$ is equal to $\widetilde r_n $ using Eq. \eqref{eqlambdaeg}.

\newcommand{\url}{\texttt}

\providecommand{\bysame}{\leavevmode ---\ }
\providecommand{\og}{``}
\providecommand{\fg}{''}
\providecommand{\smfandname}{\&}
\providecommand{\smfedsname}{\'eds.}
\providecommand{\smfedname}{\'ed.}
\providecommand{\smfmastersthesisname}{M\'emoire}
\providecommand{\smfphdthesisname}{Th\`ese}


\begin{thebibliography}{10}

\bibitem{Apery}
{\scshape R.~Ap{\'e}ry} -- {\og Irrationalit\'e de $\zeta(2)$ et
  $\zeta(3)$\fg}, in \emph{Journ{\'e}es Arithm{\'e}tiques (Luminy, 1978)},
  Ast{\'e}risque, no.~61, 1979, p.~11--13.

\bibitem{BR}
{\scshape K.~Ball {\normalfont \smfandname} T.~Rivoal} -- {\og
  Irrationalit{\'e} d'une infinit{\'e} de valeurs de la fonction z\^eta aux
  entiers impairs\fg}, \emph{Invent. Math.} \textbf{146} (2001), no.~1,
  p.~193--207.

\bibitem{DBShid}
{\scshape D.~Bertrand} -- {\og Le th\'eor\`eme de {S}iegel-{S}hidlovsky
  revisit\'e\fg}, in \emph{Number theory, Analysis and Geometry: in memory of
  Serge Lang} (\mbox{D. Goldfeld {\em et al.}}, eds.), Springer, 2012,
  p.~51--67.

\bibitem{BB}
{\scshape D.~Bertrand {\normalfont \smfandname} F.~Beukers} -- {\og
  \'{E}quations diff\'erentielles lin\'eaires et majorations de
  multiplicit\'es\fg}, \emph{Ann. Sci. \'Ecole Norm. Sup. (4)} \textbf{18}
  (1985), no.~1, p.~181--192.

\bibitem{Colmez}
{\scshape P.~Colmez} -- {\og Arithm{\'e}tique de la fonction z\^eta\fg}, in
  \emph{Journ{\'e}es math{\'e}matiques X-UPS 2002}, {\'e}ditions de l'{\'e}cole
  {P}olytechnique, 2003, p.~37--164.

\bibitem{EMS}
{\scshape N.~Fel'dman {\normalfont \smfandname} Y.~Nesterenko} -- \emph{Number
  theory {IV}, transcendental numbers}, Encyclopaedia of Mathematical Sciences,
  no.~44, Springer, 1998, A.N. Parshin and I.R. Shafarevich, eds.

\bibitem{SFdistrib}
{\scshape S.~Fischler} -- {\og Distribution of irrational zeta values\fg},
  \emph{Bull. Soc. Math. France} \textbf{145} (2017), no.~3, p.~381--409.

\bibitem{SFcaract}
\bysame , {\og Shidlovsky's multiplicity estimate and irrationality of zeta
  values\fg}, \emph{J. Austral. Math. Soc.} \textbf{105} (2018), no.~2,
  p.~145--172.

\bibitem{gfndio2}
{\scshape S.~Fischler {\normalfont \smfandname} T.~Rivoal} -- {\og Linear
  independence of values of ${G}$-functions, {II}. {O}utside the disk of
  convergence\fg}, preprint arXiv 1811.08758 [math.NT], soumis.

\bibitem{FR}
\bysame , {\og Approximants de {P}ad{\'e} et s{\'e}ries
  hyperg{\'e}om{\'e}triques {\'e}quilibr{\'e}es\fg}, \emph{J. Math. Pures
  Appl.} \textbf{82} (2003), no.~10, p.~1369--1394.

\bibitem{FSZ}
{\scshape S.~Fischler, J.~Sprang {\normalfont \smfandname} W.~Zudilin} -- {\og
  Many odd zeta values are irrational\fg}, preprint arXiv 1803.08905 [math.NT],
  Compositio Math., to appear, 2018.

\bibitem{ChowlaMilnor}
{\scshape S.~Gun, M.~R. Murty {\normalfont \smfandname} P.~Rath} -- {\og On a
  conjecture of {C}howla and {M}ilnor\fg}, \emph{Canad. J. Math.} \textbf{63}
  (2011), no.~6, p.~1328--1344.

\bibitem{HW}
{\scshape G.~Hardy {\normalfont \smfandname} E.~Wright} -- \emph{An
  introduction to the theory of numbers}, fifth \smfedname, Oxford Science
  Publications, 1979.

\bibitem{Pilehroodssums}
{\scshape T.~Hessami~Pilehrood {\normalfont \smfandname} K.~Hessami~Pilehrood}
  -- {\og Irrationality of sums of zeta values\fg}, \emph{Mat. Zametki [Math.
  Notes]} \textbf{79} (2006), no.~4, p.~607--618 [561--571].

\bibitem{Marcovecchio}
{\scshape R.~Marcovecchio} -- {\og Linear independence of linear forms in
  polylogarithms\fg}, \emph{Annali Scuola Norm. Sup. Pisa} \textbf{V} (2006),
  no.~1, p.~1--11.

\bibitem{Matala-Aho}
{\scshape T.~Matala-aho} -- {\og On {D}iophantine approximations of the
  solutions of $q$-functional equations\fg}, \emph{Proc. Roy. Soc. Edinburgh
  Sect. A} \textbf{132} (2002), p.~639--659.

\bibitem{Nash}
{\scshape M.~H. Nash} -- {\og Special values of {H}urwitz zeta functions and
  {D}irichlet ${L}$-functions\fg}, Ph.{D}. thesis, Univ. of Georgia, Athens,
  U.S.A., 2004.

\bibitem{Neukirch}
{\scshape J.~Neukirch} -- \emph{Algebraic number theory}, Springer, 1999.

\bibitem{Nishimoto}
{\scshape M.~Nishimoto} -- {\og On the linear independence of the special
  values of a {D}irichlet series with periodic coefficients\fg}, preprint arXiv
  1102.3247 [math.NT], 2011.

\bibitem{RivoalCRAS}
{\scshape T.~Rivoal} -- {\og La fonction z\^eta de {R}iemann prend une
  infinit{\'e} de valeurs irrationnelles aux entiers impairs\fg}, \emph{C. R.
  Acad. Sci. Paris, Ser. I} \textbf{331} (2000), no.~4, p.~267--270.

\bibitem{Catalan}
{\scshape T.~Rivoal {\normalfont \smfandname} W.~Zudilin} -- {\og Diophantine
  properties of numbers related to {C}atalan's constant\fg}, \emph{Math.
  Annalen} \textbf{326} (2003), no.~4, p.~705--721.

\bibitem{Shidlovski}
{\scshape A.~B. Shidlovsky} -- \emph{Transcendental numbers}, de Gruyter
  Studies in Math., no.~12, de Gruyter, Berlin, 1989.

\bibitem{Sprang}
{\scshape J.~Sprang} -- {\og Infinitely many odd zeta values are irrational.
  {B}y elementary means\fg}, preprint arXiv:1802.09410 [math.NT], 2018.

\bibitem{Zudilintrick}
{\scshape W.~Zudilin} -- {\og One of the odd zeta values from $\zeta(5)$ to
  $\zeta(25)$ is irrational. {B}y elementary means\fg}, \emph{SIGMA}
  \textbf{14} (2018), no.~028, 8 pages.

\end{thebibliography}
\end{document}